\documentclass[12pt]{article}
\usepackage[latin1]{inputenc}
\usepackage{amsmath, amssymb, eucal, mathrsfs}
\newtheorem{theorem}{Theorem}
\newtheorem{proposition}[theorem]{Proposition}
\newtheorem{lemma}[theorem]{Lemma}

\newcommand{\R}{\mathbb{R}}
\newcommand{\U}{\mathcal{U}}
\newcommand{\Le}{\mathbb{L}}

\newcommand{\C}{\mathbb{C}}

\newcommand{\spa}{\mbox{span\,}}

\newcommand{\Ima}{\mbox{Im\,}}

\newcommand{\Ve}{\mathbb{V}}

\newcommand{\po}{{\hspace*{-1ex}}{\bf .  }}

\newcommand{\nap}{\nabla^{\perp}}
\def\<{{\langle}}
\def\>{{\rangle}}
\def\lp{{\langle\!\langle}}
\def\rp{{\rangle\!\rangle}}

\def\a{\alpha}
\def\be{\begin{equation} } 
\def\ee{\end{equation} }

\def\nap{\nabla^\perp}
\def\proof{\noindent\emph{Proof: }}
\def\qed{\ifhmode\unskip\nobreak\fi\ifmmode\ifinner
\else\hskip5 pt \fi\fi\hbox{\hskip5 pt \vrule width4 pt
height6 pt  depth1.5 pt \hskip 1pt }}

\begin{document}

\title{Conformal Kaehler Euclidean submanifolds}
\author{A. de Carvalho, S. Chion and M. Dajczer}
\date{}
\maketitle

\begin{abstract} 
Let $f\colon M^{2n}\to\R^{2n+\ell}$, $n \geq 5$, denote a conformal 
immersion into Euclidean space with codimension $\ell$ of a Kaehler 
manifold of complex dimension $n$ and free of flat points.  For 
codimensions $\ell=1,2$ we show that such a submanifold can always 
be locally  obtained in a rather simple way, namely, from an isometric 
immersion of the Kaehler manifold $M^{2n}$ into either $\R^{2n+1}$ or 
$\R^{2n+2}$, the latter being a class of submanifolds already extensively 
studied.  
\end{abstract}

\renewcommand{\thefootnote}{\fnsymbol{footnote}} 
\renewcommand{\thefootnote}{\arabic{footnote}} 

\renewcommand{\thefootnote}{\fnsymbol{footnote}} 
\renewcommand{\thefootnote}{\arabic{footnote}}  

Throughout the paper $f\colon M^{2n}\to\R^{2n+\ell}$ denotes a 
\emph{conformal Kaehler submanifold}, that is, $(M^{2n},J)$ is a 
connected Kaehler manifold of complex dimension $n\geq 2$ and $f$ 
a conformal immersion into Euclidean space with codimension 
$\ell$. That the immersion is \emph{conformal} means that there 
is a positive function $\lambda \in C^\infty(M)$ such that the 
metric induced by $f$ is related to the original Kaehler metric  
by $\<\,,\,\>_f=\lambda^2\<\,,\,\>_{M^{2n}}$. The immersion is called 
a \emph{real Kaehler submanifold} if $\lambda\equiv 1$.
Our goal is to  describe, up to a conformal congruence of the ambient 
space, the local situation of the conformal Kaehler submanifold  
submanifolds if $\ell$ is at most two. Recall that two immersions 
$f,g\colon M^n\to\R^N$  are said to be \emph{conformally congruent} 
if $g=\tau\circ f$ for some conformal (Moebius) transformation 
$\tau$ of $\R^N$.
\vspace{1ex}

In our first result, by a \emph{real Kaehler hypersurface} 
we mean a real Kaehler submanifold with codimension one of
a manifold free of flat points. These 
submanifolds have been locally classified by Dajczer and Gromoll 
\cite{DG0} by means of the Gauss parametrization in terms of a
pseudoholomorphic surface in a sphere and a smooth function on 
the surface. Florit and Zheng \cite{FZ} showed that metrically 
complete real Kaehler hypersurfaces are just cylinders over a 
surface in $\R^3$.

\begin{theorem}\po\label{main1}
Any conformal immersion  $f\colon M^{2n}\to\R^{2n+1}$, $n\geq 4$, 
of a simply connected Kaehler manifold free of flat points is 
conformally congruent to a real Kaehler hypersurface.
\end{theorem}
 
For codimension two, simple examples of conformal Kaehler 
submanifolds are obtained by composing a holomorphic hypersurface 
$M^{2n}\to\C^{n+1}$ or the extrinsic product of a pair of real 
Kaehler hypersurfaces with a conformal transformation 
of the ambient space. 
But there are many other examples of real Kaehler submanifolds that 
can be composed with a conformal ambient transformation; 
for instance see Dajczer and Gromoll \cite{DG1} for the class of
complex ruled submanifolds, including the metrically complete 
that are among the ones produced by a Weierstrass type representation.  
See also Dajczer and Florit \cite{DF} for the case of submanifolds 
of rank two.

\begin{theorem}\po\label{main2}
Let $f\colon M^{2n}\to\R^{2n+2}$, $n\geq 5$, be a conformal Kaehler 
submanifold where $M^{2n}$ is free of flat points.  Then there is an 
open dense subset $M_0$ of $M^{2n}$ such that along any connected 
component $N^{2n}$ of $M_0 $  one of the following holds:
\begin{itemize}
\item[(i)]  $f|_N$ is conformally congruent to a real Kaehler 
submanifold $g\colon N^{2n}\to\R^{2n+2}$.
\item[(ii)]  $f|_N=h\circ g$ is a composition of a real Kaehler 
hypersurface $g\colon N^{2n}\to\R^{2n+1}$ and a conformal immersion 
$h\colon V\to\R^{2n+2}$ where $V\subset\R^{2n+1}$ is open and 
$g(N)\subset V$.
\end{itemize}
\end{theorem}

Notice that $h$ in part $(ii)$ is just a conformally flat hypersurface. 
The submanifolds in this class have been parametrically 
described by do Carmo, Dajczer and Mercuri \cite{CDM}.

\section{Preliminaries}

\subsection{The isometric light-cone representative}

The \emph{light-cone} $\Ve^{m+1}$ of the standard flat Lorentzian 
space $\Le^{m+2}$ is one of the two connected components of the set 
of all light-like vectors, that is,  
$$
\{v\in\Le^{m+2}:\<v,v\> =0,\;v\neq 0\} 
$$
endowed with the degenerate metric inherited from $\Le^{m+2}$. 

The Euclidean space $\R^m$ can be realized as an umbilic 
hypersurface of $\Ve^{m+1}$ as follows: 
Given light-like vectors $v,w\in \Le^{m+2}$ such that $\<v,w\>=1$ 
and a linear isometry $C\colon\R^m\to\{v,w\}^\perp$, define 
$\Psi\colon\R^m\to\Ve^{m+1}\subset\Le^{m+2}$ by 
$$
\Psi(x)=v+Cx-\frac{1}{2}\|x\|^2w.
$$
Then $\Psi$ is an isometric embedding of $\R^m$ as an umbilical 
hypersurface in the light cone given as an intersection of $\Ve^{m+1}$ 
with an affine hyperplane, namely,
$$
\Psi(\R^m)=\{y\in\Ve^{m+1}\colon\<y,w\>=1\}.
$$
The normal bundle of $\Psi$ is $N_\Psi\R^m=\spa\{\Psi,w\}$ and the
second fundamental form is 
$$
\a^\Psi(X,Y)=-\<X,Y\>_{\R^m}w.
$$
We observe that $\Psi(\R^m)$ is independent of the triple 
$v,w,C$ in the sense that different triples produce submanifolds 
congruent by an isometry of $\Le^{m+2}$.

 If $f\colon M^n\to\R^m$ is a conformal immersion with conformal 
factor $\lambda\in C^\infty(M)$, then the isometric immersion 
$$
F=\frac{1}{\lambda}\Psi\circ f\colon M^n\to\Ve^{m+1}\subset\Le^{m+2},
$$ 
is called the \emph{isometric light-cone representative} of $f$.  
The normal bundle of $F$ decomposed orthogonally as 
$$
N_FM=\Psi_*N_fM\oplus L^2
$$ 
such that $F\in\Gamma(L^2)$ and the second fundamental form 
of $F$ satisfies 
\be\label{sff}
\<\alpha^F(X,Y),F\>=-\<X,Y\>
\ee 
for all tangent vector fields $X,Y\in\mathfrak{X}(M)$.
The full expression of the second fundamental form of $F$ as well 
as additional information on the isometric light-cone representatives 
can be found in \cite{DT} and \cite{To}.

\begin{proposition}\po\label{cofcong}
Two conformal immersions $f,g\colon M^n\to\R^m$ are conformally 
congruent if and only if their isometric light-cone representatives
$F,G\colon M^n\to\Ve^{m+1}\subset\Le^{m+2}$ are isometrically 
congruent.  	
\end{proposition}

\proof See Proposition 9.18 in \cite{DT}.\qed 

\begin{proposition}\po\label{deltact}
Let $F\colon M^n\to\Ve^{n+p+1}\subset\Le^{n+p+2}$ be an isometric 
immersion that carries  a normal light-like vector field $\delta$
that is constant in $\Le^{n+p+2}$ and satisfies \mbox{$\<F,\delta\>=1$}.  
If $M^n$ is simply connected there exists an isometric immersion 
$f\colon M^n\to\R^{n+p}$ that has $F$ as its isometric light-cone 
representative.
\end{proposition}

\proof With respect to the orthogonal splitting 
$N_FM=\spa\{\delta,F\}\oplus L$ we have
\be\label{dalpha}
\alpha^F(X,Y)=-\<X,Y\>\delta+\alpha_L(X,Y)
\ee 
where $\alpha_L=\pi_L\circ\alpha^F$. 
Clearly $\alpha_L\colon TM\times TM\to L$ satisfies the Gauss equation. 
In fact, being $F$ a normal vector field parallel in the normal connection 
and $\delta$ constant, it is easy to see that $\alpha_L$ also satisfies the 
Codazzi and Ricci equations when $L$ is taken with the induced connection 
$(\nap)_L$ from $N_FM$. Hence, there are an isometric immersion 
$f\colon M^n\to\R^{n+p}$ and a vector bundle isometry 
$\phi\colon L\to N_fM$ such that 
\be\label{four}
^f\nap\circ\phi=\phi\circ(\nap)_L \;\;\text{and}\;\; 
\alpha^f=\phi\circ\alpha_L.
\ee

Let $G=\Psi\circ f\colon M^n\to\Ve^{n+p+1}\subset\Le^{n+p+2}$ 
be the isometric light-cone representative of $f$. 
Then $N_GM=\spa\{G,w\}\oplus\Psi_*N_fM$ and 
\be\label{alpha2}
\alpha^G(X,Y)=-\<X,Y\>w+\Psi_*\alpha^f(X,Y).
\ee 
Let $T\colon N_GM\to N_FM$ be the vector bundle isometry
defined by $Tw=\delta$, $T\circ G=F$ and 
$T\Psi_*\xi=\phi^{-1}\xi$ for $\xi\in N_fM$.  
Using  \eqref{dalpha}, \eqref{four} and \eqref{alpha2} 
we obtain that $T\circ\alpha^G=\alpha^F$ and 
$^F\nabla^\perp\circ T=T{}\circ^G\!\nabla^\perp$.
Hence $F$ and $G$ are isometrically congruent.\qed

\subsection{Flat bilinear forms}

Let $V^n$ and $W^{p,p}$ be real vector spaces
of dimensions $n$ and $2p$, respectively, where the latter is 
endowed with an inner product of signature $(p,p)$. 
This means that $p$ is the dimension of the subspaces of maximal 
dimension where the inner product is either positive or negative 
definite. 
A vector subspace $L\subset W^{p,p}$ is called \emph{degenerate} if 
$L\cap L^\perp\neq\{0\}$ and \emph{nondegenerate} otherwise.  

A bilinear form $\beta\colon V^n\times V^n\to W^{p,p}$ (maybe not symmetric) 
is called \emph{flat} if
$$
\<\beta(X,Y),\beta(Z,T)\>-\<\beta(X,T),\beta(Z,Y)\>=0
$$
for all $X,Y,Z,T\in V^n$. It is said that $\beta$ is \emph{null} when
$$
\<\beta(X,Y),\beta(Z,T)\>=0
$$
for all $X,Y,Z,T\in V^n$. Thus null bilinear forms are 
trivially flat. We denote
$$
\mathcal{S}(\beta)=\spa\{\beta(X,Y)\colon X,Y\in V^n\}
$$ 
and
$$
\mathcal{N}(\beta)=\{Y\in V^n\colon\beta(X,Y)=0
\;\mbox{for all}\; X\in V^n\}.
$$ 

\begin{proposition}\po\label{costum}	
Let $V^n$ and $U^p$, $2p<n$ and $1\leq p\leq 5$,  be real 
vector spaces such that there is $J\in End(V)$ satisfying 
$J^2=-I$ and $U^p$ has an inner product of any signature. 
Let $\a\colon V^n\times V^n\to U^p$ be a symmetric bilinear 
form and let $\beta\colon V^n\times V^n\to U^p\oplus U^p$ be 
the bilinear form given by 
\be\label{beta}
\beta(X,Y)=(\a(X,Y),\a(X,JY)).
\ee  
Assume that $\beta$ is flat when $W^{p,p}=U^p\oplus U^p$ is 
endowed with the inner product given by
\be\label{metric}
\lp(\xi_1,\xi_2),(\eta_1,\eta_2)\rp
=\<\xi_1,\eta_1\>_{U^p}-\<\xi_2,\eta_2\>_{U^p}.
\ee
If the subspace $\mathcal{S}(\beta)$ is nondegenerate then 
$\dim\mathcal{N}(\beta)\geq n-2p$.
\end{proposition}

\proof For $p\leq 5$ the proof of Proposition $10$ in \cite{CCD}, 
where the inner product on $U^p$ is positive definite and $\alpha$ 
satisfies a certain condition, can be adapted to this case. With 
the notations in there and using the same  type of arguments used
there it is easy to conclude that the only cases one needs to 
consider are $2\leq\kappa<\tau\leq p-1$
where $\kappa$ and $\tau$ are even.  Thus, we only have to deal 
with the case $\tau=4$ and $\kappa=2$.  

By Fact $11$ in \cite{CCD} there exist 
$Y_1,Y_2\in RE^o(\beta)\cap RE(\hat{\beta})$ such that
$$
\hat{\U}(X)=\hat{B}_{Y_1}(V)+\hat{B}_{Y_2}(V).
$$
We cannot have $B_{Y_j}(N)=\U(X)$, $j=1,2$, since otherwise 
$\dim U(Y_j)\leq 3<\tau$.  Therefore, by Fact $12$ in \cite{CCD} it 
remains to consider the case $\dim B_{Y_j}(N)\leq 2$, $j=1,2$.  Set 
$B_1=B_{Y_1}|_N\colon N\to\U(X)$, $N_1=\ker B_1$, 
$B_2=B_{Y_2}|_{N_1}\colon N_1\to\U(X)$ and 
$N_2=\ker B_2$.  Then $N_2\subset\mathcal{N}(\beta)$ and 
$$
\dim\mathcal{N}(\beta)\geq\dim N_2\geq\dim N_1-2\geq\dim N-4\geq n-2p,
$$
and this concludes the proof.\qed

\section{The proofs}

The following application of Proposition \ref{costum} is the 
main ingredient in the proofs of the theorems in this paper.

\begin{proposition}\po\label{main}
Let $V^n$ and $U^p$, $2p<n$ and $1\leq p\leq 5$,  be real vector 
spaces such that there is $J\in End(V)$ satisfying $J^2=-I$ and 
$U^p$ carries an either positive definite or Lorentzian inner product.
Assume that the bilinear form 
$\beta\colon V^n\times V^n\to W^{p,p}=U^p\oplus U^p$ defined by 
\eqref{beta} is flat with respect to the inner product 
\eqref{metric}. If $\dim\mathcal{N}(\beta)\leq n-2p-1$  
then $\U=\mathcal{S}(\beta)\cap\mathcal{S}(\beta)^\perp$ satisfies 
$\dim\U=s>0$ is even. Moreover, let $L\subset U^p$ denote the 
projection of $\U$ on the first factor of $W^{p,p}$. Then, 
we have:
\begin{itemize}
\item[(i)] If the subspace $L$ is nondegenerate then $\dim L=s$ 
and $L$ inherits a positive definite inner product. With respect to 
the orthogonal splitting  $U^p=L\oplus L^\perp$ we denote 
$\a_1=\pi_L\circ\a$ and $\a_2=\pi_{L^\perp}\circ\a$. Then 
$$
\a_1(X,JY)=\a_1(JX,Y)\;\;\mbox{for all}\;\;X,Y\in V^n
$$
and 
$$
\dim\mathcal{N}(\a_2)\cap J\mathcal{N}(\a_2)\geq n-2(p-s).
$$
\item[(ii)] If the subspace $L$ is degenerate let 
$0\neq\delta\in L\cap L^\perp$.  Then there is an orthogonal 
splitting $U^p=U_0\oplus U_1^{s-2}\oplus U_2^{p-s}$, 
$s=2\;\text{or}\;4$,  with $U_0=\spa\{\delta,\zeta\}$, where 
$\zeta\in U^p$ is a light-like vector satisfying 
$\<\delta,\zeta\>=1$,  and 
$L=\spa\{\delta\}\oplus U_1^{s-2}$ such that
$\a_j=\pi_{U_j}\circ\a$, $j=0,1,2$, satisfy
$$
\<\a(X,Y),\delta\>=0\;\;\text{and}\;\;
\a_1(X,JY)=\a_1(JX,Y)\;\;\mbox{for all}\;\;X,Y\in V^n
$$
and
$$
\dim\mathcal{N}(\a_2)\cap J\mathcal{N}(\a_2)\geq n-2(p-s).
$$ 
\end{itemize}
\end{proposition}

\proof By Proposition \ref{costum} we have $s>0$. 
If $0\neq(\xi,\bar{\xi})\in\U$, then
$$
(\xi,\bar{\xi})=\sum_i\beta(X_i,Y_i)
=\sum_i(\a(X_i,Y_i),\a(X_i,JY_i))
$$
and
$$
0=\lp\beta(X,Y),(\xi,\bar{\xi})\rp=\<\a(X,Y),\xi\>
-\<\a(X,JY),\bar{\xi}\>
$$
for any $X,Y\in V^n$.  Then 
$(\bar{\xi},-\xi)=\sum_i\beta(X_i,JY_i)\in\mathcal{S}(\beta)$
and 
$$
\lp(\beta(X,Y),(\bar{\xi},-\xi)\rp=\<\a(X,Y),\bar{\xi}\>
+\<\a(X,JY),\xi\>=0
$$
for any $X,Y\in V^n$. Thus also $(\bar{\xi},-\xi)\in\U$. 
It follows that $s$ is even and  
$$
\pi_1(\U)=L=\pi_2(\U),
$$ 
where $\pi_j\colon W^{p,p}\to U^p$, $j=1,2$, 
denote the projections onto the factors.
\vspace{1ex}

\noindent {\emph Case $(i)$}: We have that the inner product induced 
on $L$ is positive definite. In fact, if otherwise there are vectors 
$\delta,\bar{\delta}\in L$ such that $\delta$ is time-like and 
$(\delta,\bar{\delta}), (\bar{\delta},-\delta)\in\U$.
But then also $\bar{\delta}$ would be a time-like vector orthogonal 
to $\delta$ in contradiction with the signature of $U^p$.

We have $\beta=\beta_1+\beta_2$ where 
$$
\beta_j(X,Y)=(\a_j(X,Y),\a_j(X,JY)),\;j=1,2.
$$
Since $\beta_1$ is null, then
$$
0=\lp\beta_1(X,Y),\beta_1(Z,W)\rp
=\<\a_1(X,Y),\a_1(Z,W)\>-\<\a_1(X,JY),\a_1(Z,JW)\>.
$$
Then   
$T\colon\mathcal{S}(\a_1)\to\mathcal{S}(\a_1)$ defined by
$$
T\a_1(X,Y)=\a_1(X,JY)
$$  
is a linear isometry and
$$
\a_1(JX,Y)=\a_1(Y,JX)=T\a_1(Y,X)
=T\a_1(X,Y)=\a_1(X,JY).
$$
Being $\beta$ flat and $\beta_1$ null, then also $\beta_2$ 
is flat. Since the subspace $\mathcal{S}(\beta_2)$ 
is nondegenerate we have from Proposition \ref{costum} that
$$
\dim\mathcal{N}(\beta_2)\geq n-2(p-s).
$$
To conclude the proof of this case observe that
$\mathcal{N}(\beta_2)=\mathcal{N}(\a_2)\cap J\mathcal{N}(\a_2)$.
\vspace{1ex}

\noindent {\emph Case $(ii)$}:
Let $\bar{\delta}\in L$ be such that 
$(\delta,\bar{\delta}),(\bar{\delta},-\delta)\in\U$.  
Since the inner product on $U^p$ has Lorentzian signature,  
then the vectors $\delta,\bar{\delta}$ must be linearly 
dependent. Thus $(\delta,0),(0,\delta)\in\U$, and hence 
$$
0=\lp\beta(X,Y),(\delta,0)\rp=\<\a(X,Y),\delta\>
$$
for all $X,Y\in V^n$.

Assume $s=2$, in which case $\U=\spa\{(\delta,0),(0,\delta)\}$. 
We have $\a=\a_0+\a_2$, where
$$
\a_0(X,Y)=\<\a(X,Y),\zeta\>\delta.
$$
Hence  $\beta={\beta}_0+\beta_2$, where 
$$
\beta_0(X,Y)=\<\a(X,Y),\zeta\>(\delta,0)+\<\a(X,JY),\zeta\>(0,\delta)
$$
and
$$
\beta_2(X,Y)=(\a_2(X,Y),\a_2(X,JY)).
$$  
Because $\beta$ is flat and $\beta_0$ is null, then
$\beta_2$  is flat. Since the subspace $\mathcal{S}(\beta_2)$ is 
nondegenerate, then Proposition \ref{costum} gives
$$
\dim\mathcal{N}(\beta_2)\geq n-2p+4.
$$

Assume $s=4$. Then there are space-like vectors 
$\xi,\bar{\xi}\in L$ such that
$$
\U=\spa\{(\delta,0),(0,-\delta),(\xi,\bar\xi),(\bar\xi,-\xi)\}.
$$
Set $U_1=\spa\{\xi,\bar{\xi}\}$ and choose $\zeta\perp U_1$.
Let $\beta_j\colon V^n\times V^n\to U_j\oplus U_j$ be given by 
$$
\beta_j(X,Y)=(\a_j(X,Y),\a_j(X,JY)),\;j=0,1,2.
$$
Then $\beta=\beta_0+\beta_1+\beta_2$ where $\beta_0$ and 
$\beta_1$ are null.

If $T\colon U_1\to U_1$ be the linear isometry defined by 
$$
T\a_1(X,Y)=\a_1(X,JY),
$$
then
$$
\a_1(JX,Y)=\a_1(Y,JX)=T\a_1(Y,X)
=T\a_1(X,Y)=\a_1(X,JY).
$$
Since  $\beta_2$ is flat and $\mathcal{S}(\beta_2)$ is  a
nondegenerate subspace, then Proposition \ref{costum} gives
$$
\dim\mathcal{N}(\beta_2)\geq n-2p+8,
$$
and this concludes the proof.\vspace{2ex}\qed
 
 Let $f\colon M^{2n}\to\R^{2n+p}$ be a conformal immersion of 
a simply connected Kaehler manifold free of flat points and 
let $\a^F\colon TM\times TM\to N_FM$ be the second fundamental 
form of its isometric light-cone representative 
$F\colon M^{2n}\to\Ve^{2n+p+1}\subset\Le^{2n+p+2}$. 
At any point $x\in M^{2n}$, let 
$\beta\colon T_xM\times T_xM\to W^{p+2,p+2}=N_FM(x)\oplus N_FM(x)$ 
be the bilinear form given by 
\be\label{betak}
\beta(X,Y)=(\a^F(X,Y),\a^F(X,JY))
\ee  
where the inner product in $W^{p+2,p+2}$ is as in \eqref{metric}.
Using the Gauss equation and that the curvature tensor of
$M^{2n}$ satisfies $J\circ R(X,Y)=R(X,Y)\circ J$
it is easy to verify that $\beta$ is flat. Moreover,  since  
$$
\<\a^F(X,Y),F \>=-\<X,Y\>
$$
we have $\mathcal{N}(\beta)\subset\mathcal{N}(\a^F)=\{0\}$. 
\vspace{2ex}

\noindent \emph{Proof of Theorem \ref{main1}:} From 
Proposition \ref{main} applied at $x\in M^{2n}$ to  
$\beta\colon T_xM\times T_xM\to W^{3,3}$ defined  
by \eqref{betak} in terms of $\a^F$ satisfying \eqref{sff} 
we have $s(x)=2$. We also have that  $L$ is a degenerate subspace.  
In fact, if otherwise, by part $(i)$ there exists an orthogonal 
splitting $N_FM(x)=L\oplus L^\perp$ such that $L$ inherits a positive 
definite inner product and $L^\perp=\spa\{\eta\}$ where 
$\eta\in N_FM(x)$ is a unit time-like vector. 
Moreover, the $J$-invariant subspace 
$\Delta=\mathcal{N}(\beta_2)$ satisfies $\dim\Delta\geq 2n-2$
and, since $L=\spa\{\xi,\bar\xi\}$ where $(\xi,\bar\xi)\in\U$, 
then the shape operators of $F$ satisfy
$$
A_\xi^F=-JA_{\bar\xi}^F\;\;\mbox{and}\;\;\Delta\subset\ker A_\eta^F.
$$ 
Since $F\in N_FM(x)$ then $F=a\xi+b\bar\xi+c\eta$. Hence
$$
JZ=-JA_F^FZ=-J(aA_\xi^FZ+bA_{\bar{\xi}}^FZ)
=A^F_{b\xi-a\bar\xi}Z
$$
for any $Z\in\Delta$. Therefore $J|_\Delta
=A^F_{b\xi-a\bar\xi}|_\Delta$, and this is a contradiction.

Since the subspace $L$ is degenerate, by part $(ii)$ at any 
point there is a splitting
$$
N_FM=\spa\{\delta,\zeta\}\oplus U_2
$$
such that $A_\delta^F=0$ and the $J$-invariant subspace 
$\Delta=\mathcal{N}(\beta_2)$ satisfies $\dim\Delta\geq 2n-2$. 
Moreover, since $M^{2n}$ is free of flat points then 
$\dim\Delta=2n-2$.  

Because $F,\delta\in N_FM(x)$ are linearly independent we may
take $\zeta=F$. Hence, we have a normal basis $\{\delta,F,\xi\}$ 
with $\xi\perp\spa\{\delta,F\}$ of unit length such that
$$
A_\delta^F=0,\;\; A_F^F=-I\;\;\text{and}\;\;\Delta\subset\ker A_\xi^F.
$$
We have that $\U=\mathcal{S}(\beta)\cap\mathcal{S}(\beta)^\perp$ 
has constant dimension and hence is smooth. It follows easily that
also the frame $\{\delta,F,\xi\}$ can be taken to be smooth.

The Codazzi equation for 
$A_\delta^F$ is
$$
A_{\nabla_X^\perp \delta}^FY=A_{\nabla_Y^\perp\delta}^FX
$$
for all $X,Y\in TM$. Using that $M^{2n}$ does not have flat points, 
it is not difficult to conclude that $\delta$ is parallel in the 
normal connection, and hence constant in the ambient space. 
Thus, as in the proof of Proposition \ref{deltact}, 
there  exists an isometric immersion $g\colon M^{2n}\to\R^{2n+1}$ 
that has $A^g=A_\xi^F$ as shape operator
and its isometric light-cone representative $G=\Psi\circ g$ is 
isometrically congruent to $F$.  Therefore, by 
Proposition~\ref{cofcong}, $f$ and $g$ are conformal.
\vspace{2ex}\qed

For the proof of Theorem \ref{main2} we need the following
two technical results.

\begin{lemma}\po\label{confcompg}
Let $g\colon M^n\to\R^{n+1}$ be an isometric immersion and let 
$f\colon M^n\to\R^{n+p}$ be a conformal immersion. Then let 
$G\colon M^n\to\Ve^{n+2}\subset\Le^{n+3}$ and 
$F\colon M^n\to\Ve^{n+p+1}\subset\Le^{n+p+2}$ be the isometric
light-cone representatives of $g$ and $f$, respectively. 
Given an open subset $U\subset M^n$, there exists a conformal 
immersion $h\colon V\to\R^{n+p}$ of an open subset 
$V\supset g(U)$ of $\R^{n+1}$ such that $f|_U=h\circ g|_U$ 
if and only if there exists an isometric immersion 
$H\colon W\to\Ve^{n+p+1}$ of an open subset 
$W\subset\Ve^{n+2}$ with $G(U)\subset W$
such that $F|_U=H\circ G|_U$.
\end{lemma}

\proof See Proposition 2 in \cite{DT2}.\qed 

\begin{lemma}\po\label{le:blem2} 
Let $F\colon M^n\to\Ve^{n+3}\subset\Le^{n+4}$ be an isometric 
immersion and let $\xi$ be a normal vector field of unit length 
that satisfies $\<\xi,F\>=0$, $\text{rank}\,A_\xi^F=1$ and 
is parallel along $\ker A_\xi^F$. Then, there exist open 
subsets $V\subset M^n$ and $W\subset\Ve^{n+2}$  and local
isometric immersions $G\colon V\to\Ve^{n+2}$ and 
$H\colon W\to\Ve^{n+3}$ with $G(V)\subset W$ 
such that $F|_V=H\circ G$.
\end{lemma}

\proof See Lemma 2 in \cite{CT}.\vspace{2ex}\qed

\noindent\emph{Proof of Theorem \ref{main2}:} 
We proceed making use of the definitions and notations in 
the proof of Theorem \ref{main1}. Proposition \ref{main} 
applied to $\beta\colon T_xM\times T_xM\to W^{4,4}$ at 
$x\in M^{2n}$ gives $s(x)=2\,\text{or}\,4$. 
In what follows we work on an open dense subset $M_*$ of 
$M^{2n}$ where $\dim\mathcal{S}(\beta)$ is locally constant.  
Let $M_2$ be the open subset of $M_*$ where $s(x)=2$.
A similar argument as in the proof of Theorem \ref{main1} gives 
that the subspace $L$ is degenerate at each point of $M_2$.  
Then along $M_2$ there is a smooth orthogonal splitting of 
the normal bundle of $F$ as
$$
N_FM=\spa\{\delta,F\}\oplus P
$$
with $\<\delta,F\>=1$ such that $A_\delta^F=0$, $A_F^F=-I$ and $\Delta
=\mathcal{N}(\beta_2)$ satisfies $\dim\Delta\geq 2n-4$.

Let $P=\spa\{\xi_1,\xi_2\}$ where the smooth frame is orthonormal.
In the sequel, we work on a connected component $M_2'$ of the open 
subset of $M_2$ where $\dim\Delta$ and the ranks of the $A_{\xi_j}^F$'s 
are locally constant.  The corresponding Codazzi equation are 
\begin{align}\label{confcod2}
\nabla_XA^F_{\xi_i}Y
&-A^F_{\xi_i}\nabla_XY
+\<\nabla_X^\perp\xi_i,\delta\>Y
-\<\nabla_X^\perp\xi_i,\xi_j\>A_{\xi_j}^F Y\\
&=\nabla_YA_{\xi_i}^FX-A_{\xi_i}^F\nabla_YX
+\<\nabla_Y^\perp\xi_i,\delta\>X
-\<\nabla_Y^\perp\xi_i,\xi_j\>A_{\xi_j}^F X,
\;1\leq i\neq j\leq 2,\nonumber
\end{align}
for any $X,Y\in\mathfrak{X}(M)$.  It follows that
$$
A_{\xi_j}^F[S,T]=\<\nabla_S^\perp\xi_j,\delta\>T
-\<\nabla_T^\perp\xi_j,\delta\>S
$$
for any $S,T\in\Gamma(\Delta)$.  Hence
$\<\nabla_S^\perp\xi_j,\delta\>=0$, $j=1,2$ and 
$S\in\Gamma(\Delta)$. Then \eqref{confcod2} yields
$$
-A_{\xi_i}^F\nabla_XS+\<\nabla_X^\perp{\xi_i},\delta\>S
=\nabla_SA_{\xi_i}^F X-A^F_{\xi_i}\nabla_SX
-\<\nabla_S^\perp\xi_i,\xi_j\>A^F_{\xi_i}X,\; 1\leq i\neq j\leq 2,
$$
for any $S\in\Gamma(\Delta)$ and  $X\in\mathfrak{X}(M)$. 
In particular, 
\be\label{2umb}
\<\nabla_X^\perp\delta,\xi_j\>\<S,T\>
=\<\nabla_ST,A_{\xi_j}^F X\>,\;j=1,2,
\ee
for any $X\in\mathfrak{X}(M)$ and $S,T\in\Gamma(\Delta)$.
Thus 
\be\label{cumb3}
(\nabla_ST)_{\Ima A_{\xi_j}^F}=\<S,T\>Z_j,\; j=1,2,
\ee 
where $Z_j\in\Gamma(\Ima A_{\xi_j}^F)$
and $S,T\in\Gamma(\Delta)$. Now \eqref{2umb} reads as
\be\label{econf4}
\<\nabla_X^\perp\delta,\xi_j\>=\<A_{\xi_j}^FZ_j,X\>,\; j=1,2,
\ee
for any $X\in\mathfrak{X}(M)$.

We denote
$$
R(x)=\spa\{A_{\xi_j}^FZ_j(x):x\in M_2',\; j=1,2\}.
$$
We claim that the open subset $N_2\subset M_2'$ defined by
$$
N_2=\{x\in M_2':\dim R(x)=2\}
$$
is empty. In fact, the Codazzi equation for $A_\delta^F$ is
\be\label{coddelta}
\<\nabla_X^\perp\delta,\xi_1\>A_{\xi_1}^FY
+\<\nabla_X^\perp\delta,\xi_2\>A_{\xi_2}^FY
=\<\nabla_Y^\perp\delta,\xi_1\>A_{\xi_1}^FX
+\<\nabla_Y^\perp\delta,\xi_2\>A_{\xi_2}^FX
\ee
for any $X,Y\in\mathfrak{X}(M)$.
We have from \eqref{econf4} that $\delta$ is parallel along $
R^\perp$. Hence
\be\label{reduced}
\<\nabla_Y^\perp\delta,\xi_1\>A_{\xi_1}^FX
+\<\nabla_Y^\perp\delta,\xi_2\>A_{\xi_2}^FX=0
\ee
for any $X\in\Gamma(R^\perp)$ and $Y\in\mathfrak{X}(M)$.
In particular, the vectors $A_{\xi_1}^FX, A_{\xi_2}^FX$ cannot 
be linearly independent for any $X\in R^\perp$.
If otherwise \eqref{reduced} yields that $\delta$ is parallel,
and then \eqref{econf4} gives $R=0$, a contradiction. 

We argue that 
\be\label{cont}
R^\perp\subset\ker A_{\xi_1}^F\cap\ker A_{\xi_2}^F.
\ee
Suppose that $X\notin\ker A_{\xi_1}^F\cap\ker A_{\xi_2}^F$
for $X\in\Gamma(R^\perp)$.  By the above
$A_{\xi_1}^FX=\gamma A_{\xi_2}^FX$ where $A_{\xi_2}^FX\neq 0$
but $\gamma\in C^\infty(N_2)$ may vanish. We obtain from 
\eqref{reduced} that
$$
\gamma\<\nabla_Y^\perp\delta,\xi_1\>+\<\nabla_Y^\perp\delta,\xi_2\>=0
$$
for any $Y\in\mathfrak{X}(M)$.  
Then \eqref{econf4} gives
$\gamma A_{\xi_1}^FZ_1+A_{\xi_2}^FZ_2=0$, which is 
a contradiction. 

Since $M^{2n}$ is free of flat points, we have from \eqref{cont}
that we may choose the frame $\xi_1,\xi_2$ such that 
$\ker A_{\xi_1}^F=R^\perp=\ker A_{\xi_2}^F$. From 
\eqref{cumb3} we obtain $Z_1=Z_2=Z\in R$.  
Then \eqref{econf4} and \eqref{coddelta} give 
\be\label{curve}
\<A_{\xi_1}^FZ,X\>A_{\xi_1}^FY+\<A_{\xi_2}^FZ,X\>A_{\xi_2}^FY
=\<A_{\xi_1}^FZ,Y\>A_{\xi_1}^FX+\<A_{\xi_2}^FZ,Y\>A_{\xi_2}^FX
\ee
for any $X,Y\in\mathfrak{X}(M)$. By the Gauss equation 
\eqref{curve} is equivalent to $R(X,Y)Z=0$, and this is a 
contradiction since $M^{2n}$ is free of flat points.  
Thus the claim that $N_2$ is empty has been proved.

Let $N_1\subset M_2'$ be the open subset defined as
$$
N_1=\{x\in M_2':\dim R(x)=1\}.
$$
Similarly as above, we obtain that $\delta$ is parallel along the 
hyperplane $R^\perp$ and that the vectors $A_{\xi_1}^FX, A_{\xi_2}^FX$ 
cannot be linearly independent for  any $X\in R^\perp$.  And from
\eqref{econf4} we have that $\delta$ is not a parallel vector field.
Hence, we can choose the frame $\{\xi_1,\xi_2\}$ for $P$ such that 
\be\label{null}
\<\nap_X\delta,\xi_2\>=0
\ee
for any $X\in\mathfrak{X}(M)$.
It now follows from \eqref{coddelta} that 
$R^\perp\subset\ker A_{\xi_1}^F$.

We have $R^\perp=\ker A_{\xi_1}^F$. In fact, otherwise
$A_{\xi_1}^F=0$ and the Codazzi equation gives
$$
\<\nabla_X^\perp\xi_1,\delta\>S
=-\<\nabla_S^\perp\xi_1,\xi_2\>A_{\xi_2}^FX
$$
for any $S\in\Gamma(\Delta)$ and $X\in\mathfrak{X}(M)$.  
It follows that $\delta$ is a parallel vector field, and this is 
a contradiction.

We obtain from \eqref{econf4} and \eqref{null} that $A_{\xi_2}^FZ_2=0$.  
Since $A_{\xi_2}^F|_{\Ima A_{\xi_2}^F}$ is an isomorphism, then 
$Z_2=0$.  If $\Ima A_{\xi_1}^F\subset\Ima A_{\xi_2}^F$, it follows
from \eqref{cumb3} that $Z_1=0$.  
Thus $\Ima A_{\xi_1}^F\not\subset\Ima A_{\xi_2}^F$ 
and since $M^{2n}$ has no flat points we have that
$2\leq\dim\Ima A_{\xi_2}^F\leq 3$. In fact, it holds that 
$\dim\Ima A_{\xi_2}^F=2$ since,  otherwise,
$\Delta^\perp=\Ima A_{\xi_1}^F\oplus\Ima A_{\xi_2}^F$.
Then \eqref{cumb3} gives
$$
(\nabla_ST)_{\Delta^\perp}=\<S,T\>Z_1
$$
for any $S,T\in\Delta$, and hence
$$
\<S,T\>JZ_1=J(\nabla_ST)_{\Delta^\perp}=(J\nabla_ST)_{\Delta^\perp}
=(\nabla_SJT)_{\Delta^\perp}=\<S,JT\>Z_1.
$$
Thus $Z_1=0$, and this is a contradiction.  

If $Y,Z\in\Gamma(\ker A_{\xi_1}^F)$ are linearly independent,
then the  Codazzi equation for $A_{\xi_1}^F$ is 
$$
A_{\xi_1}^F[Z,Y]
=A_{\xi_2}^F(\<\nabla_Y^\perp\xi_1,\xi_2\>Z
-\<\nabla^\perp_Z\xi_1,\xi_2\>Y).
$$
Since $\Ima A_{\xi_1}^F\not\subset\Ima A_{\xi_2}^F$ and
$A_{\xi_2}^F|_{\Ima A_{\xi_2}^F}$ is an isomorphism, then
$\<\nabla_Y^\perp\xi_1,\xi_2\>=0$ for any $Y\in\Gamma(\ker A_{\xi_1}^F)$.  
Thus $\xi_1$ is parallel along $\ker A_{\xi_1}^F$.

By Lemma \ref{le:blem2} there exist an simply connected open neighborhood 
$V\subset N_1$ of any $x\in N_1$, an open subset $W\subset\Ve^{2n+2}$ 
and local isometric immersions $G\colon V\to\Ve^{2n+2}$ 
and $H\colon W\to\Ve^{2n+3}$ with $G(V)\subset W$ such 
that $F|_V=H\circ G$. An elementary argument gives that there
exists a conformal immersion $g\colon V\to\R^{2n+1}$ that 
has $G$ as isometric light-cone representative; see pg.\ $7$ of 
\cite{To} or Proposition $9.9$ of \cite{DT}.
By Lemma \ref{confcompg} there exist a conformal immersion 
$h\colon U\to\R^{2n+2}$ such that $f|_V=h\circ g|_V$ with 
$g(V)\subset U$. Finally, by Theorem \ref{main1} applied to $g$ 
we are as in part $(ii)$ of Theorem \ref{main2}.

Let $N_0'$ be an open simply connected subset of the set 
$$
N_0=\text{int}\{x\in M_2:R(x)=0\}.
$$
Then $\delta$ is constant on $N_0'$ from \eqref{econf4}. 
By Proposition \ref{deltact}, there is an isometric 
immersion $g_0\colon N_0'\to\R^{2n+2}$ whose isometric 
light-cone representative is $F|_{N_0'}$.  From 
Proposition \ref{cofcong} we are in part $(i)$ of 
Theorem \ref{main2}.
  
To conclude, let $M_4\subset M^{2n}$ be the interior of 
the set $\{x\in M^{2n}:s(x)=4\}$.  
Then $L(x)$ for $x\in M_4$ is a degenerate subspace since, 
otherwise, $N_FM_4(x)=L(x)$ which contradicts the fact that $L(x)$ 
has a positive definite inner product. By Proposition \ref{main}, 
there is a smooth orthogonal vector bundle decomposition
$$
N_FM_4=\spa\{\delta,F\}\oplus U_1^2
$$
such that $A_\delta^F=0$, $A_F^F=-I$ and $\<\delta,F\>=1$.   
Moreover, we have $U_1^2 =\spa\{\xi_1,\xi_2\}$ where the 
smooth frame is orthonormal and
\be\label{2u4}
A_{\xi_1}^F=JA_{\xi_2}^F.
\ee
Comparing the Codazzi equations for $A_{\xi_1}^F$ and $A_{\xi_2}^F$ 
by means of \eqref{2u4}, it follows easily that 
$\delta\in\Gamma(N_FM_4)$ is parallel, hence constant in $\Le^{n+4}$. 
Along any simply connected open subset of $M_4$ we now combine 
Proposition \ref{deltact} and Lemma \ref{confcompg} to conclude 
that we are again in part $(i)$ of Theorem \ref{main2}.
Notice that in this case $f$ is minimal.\qed

\noindent Alcides de Carvalho, Sergio Chion, Marcos Dajczer\\
IMPA -- Estrada Dona Castorina, 110\\
22460--320, Rio de Janeiro -- Brazil\\
e-mail: alcidesj@impa.br, sergio.chion@impa.br, marcos@impa.br
\end{document}